\documentclass[11pt]{article}
\usepackage{amsfonts}
\usepackage{amssymb}

\newtheorem{theorem}{Theorem}[section]

\newtheorem{lemma}[theorem]{Lemma}
\newtheorem{corollary}[theorem]{Corollary}

\newtheorem{remark}{Remark}[section]
\newtheorem{example}{Example}[section]
 \usepackage{amssymb}
\usepackage{epsfig}
\usepackage{amsmath}
\usepackage{amsfonts}

\newcommand{\R}{{{\Bbb R}}}

\newcommand{\N}{{{\Bbb N}}}
\newcommand{\Nn}{{\scriptsize {\Bbb N}}}

\def\qed{\hbox to 0pt{}\hfill$\rlap{$\sqcap$}\sqcup$\medbreak}

\title{Peano's Existence Theorem revisited}
\date{February, 2012}
\begin{document}
\maketitle

\vspace{1cm}

\begin{center}
{\large Rodrigo L\'opez Pouso \\
Departamento de
An\'alise Matem\'atica\\
Facultade de Matem\'aticas,\\Universidade de Santiago de Compostela, Campus Sur\\
15782 Santiago de
Compostela, Spain.
}
\end{center}

\begin{abstract}
We present new proofs to four versions of Peano's Existence Theorem for ordinary differential equations and systems. We hope to have gained readability with respect to other usual proofs. We also intend to highlight some ideas due to Peano which are still being used today but in specialized contexts: it appears that the lower and upper solutions method has one of its oldest roots in Peano's paper of 1886.
\end{abstract}

\noindent
{\it ``Le dimostrazioni finora date dell'esistenza degli integrali delle equazioni differenziali lasciano a desiderare sotto l'aspetto della semplicit$\grave{a}$."} 

\medbreak

\noindent
{\bf G. Peano}, {\it Sull'integrabilit\`a delle equazioni differenziali di primo
ordine}, {Atti. Accad. Sci. Torino, {vol. 21} (1886), 677--685.} 

 \section{Introduction}
The fundamental importance of Peano's Existence Theorem hardly needs justification when we find it in almost any undergraduate course on the subject. Let us simply point out that Peano's Theorem provides us with a very easily checkable condition to ensure the existence of solutions for complicated systems of ordinary differential equations.
 
This paper contains a new proof to Peano's Existence Theorem and some other new proofs to not so well--known finer versions of it in the scalar case. 

\bigbreak

Before going into detail we need to introduce some notation. Let $t_0 \in \R$ and $y_0=(y_{0,1},y_{0,2},\dots, y_{0,n}) \in \R^n$ ($n \in \N$) be fixed, and let 
$$f:\mbox{Dom}(f) \subset \R^{n+1} \longrightarrow \R^n$$ be defined
in a neighborhood of $(t_0,y_0)$. We consider the initial value pro\-blem 
\begin{equation}
\label{ivp}
y'=f(t,y), \quad y(t_0)=y_0,
\end{equation}
which is really a system of $n$ coupled ordinary differential equations along with their corresponding initial conditions. 

By a solution of (\ref{ivp}) we mean a function $\varphi:I \longrightarrow \R^n$ which is diffe\-ren\-tiable in a nondegenerate real interval $I$, $t_0 \in I$, $\varphi(t_0)=y_0$, and for all $t \in I$ we have $\varphi'(t)=f(t,\varphi(t))$, the corresponding one--sided derivatives being considered at the endpoints of the interval. It is implictly required that solutions have their graphs inside the domain of the function $f$.
\bigbreak

 The most popular version of Peano's Theorem reads as follows:

\bigbreak

\noindent
{\bf Peano's Theorem.} {\it If the function $f$ is continuous in a neighborhood of $(t_0,y_0)$ then the initial value problem (\ref{ivp}) has at least one solution defined in a neighborhood of $t_0$.}

\medbreak
Peano proved the result in dimension $n=1$ first, see \cite{peano}, and then he extended it to systems in \cite{peano2}. Peano's Theorem has attracted the attention of many mathematicians and, as a result, there are many different proofs available today. We can clasify them into two fundamental types:
\begin{enumerate}
\item[(A)] Proofs based on the construction of a sequence of approximate solutions (mainly Euler--Cauchy polygons or Tonelli sequences) which converges to some solution.
\item[(B)] Proofs based on fixed point theorems (mainly Schauder's Theorem) applied to the equivalent integral version of (\ref{ivp}).
\end{enumerate}

Both types of proofs (A) and (B) have their advantages and their drawbacks. First, proofs of type (A) are more elementary and hence {\it a priori} more adequate for elementary courses. Second, proofs of type (B) are shorter and much clearer, but they involve more sophisticated results, namely fixed point theorems in infinite dimensional function spaces.

One can find some other types of proofs in the literature, but in our opinion they are not so clear or elementary. For instance, a well--known kind of mix of types (A) and (B) consists in approximating $f$ by polinomials $f_n$ (by virtue of the Stone--Weierstrass Theorem) so that the problems (\ref{ivp}) with $f$ replaced by $f_n$ have unique solutions $y_n$ (by virtue of the Banach contraction principle) which form an approximating sequence. More elementary proofs were devised in the seventies of the 20th century, partly as a reaction to a question posed in \cite{ken}. In \cite{dv, gar, jwal, wwal} we find proofs of type (A) in dimension $n=1$ which avoid the Arzel\`a--Ascoli Theorem.  As already stated in those references, similar ideas do not work in higher dimension, at least without further assumptions. We also find in \cite{dv} a proof which uses Perron's method, see \cite{perron}, a refined version of Peano's own proof in \cite{peano}. It is fair to acknowledge that probably the best application of Perron's method to (\ref{ivp}) in dimension $n=1$ is due to Goodman \cite{goo}, who even allowed $f$ to be discontinuous with respect to the independent variable and whose approach has proven efficient in more general settings, see \cite{bp, bs, has, po2}.

\bigbreak

In this paper we present a proof of Peano's Theorem which, in our opinion, takes profit from the most advantageous ingredients of proofs of types (A) and (B) to produce a new one which we find more readable. In particular, our proof in Section 2 meets the following objectives:
\begin{enumerate}
\item It involves an approximating sequence, but we do not have to worry about its convergence or its subsequences at all.
\item It involves a mapping from a space of functions into the reals, thus introducing some elements of functional analysis, but the most sophisticated result we use to study it is the Arzel\`a--Ascoli Theorem.  
\item Compactness in the space of continuous functions is conveniently emphasized as a basic ingredient in the proof, but we think that the way we use it (``continuous mappings in compact sets have minima") leads to a more readable proof than the usual (``sequences in compact sets have convergent subsequences").
\end{enumerate}

This paper is not limited to a new proof of Peano's Theorem. With a little extra work and a couple of new ideas, we prove the existence of the least and the greatest solutions for scalar problems (Section 3) and we also study the existence of solutions between given lower and upper solutions (Section 4). Comparison with the literature is discussed in relevant places and examples (some of them new) are given to illustrate the results or their limitations.

\section{Proof of Peano's Theorem}
{\bf A simplification of the problem.} We claim that {\it  it suffices to prove the existence of solutions defined on the right of the initial time $t_0$}. 

To justify it, assume we have already proven the result for that specific type of solutions and consider the problem with reversed time
\begin{equation}
\label{ivpdos}
y'=g(t,y)=-f(-t,y), \quad y(-t_0)=y_0.
\end{equation}
The function $g$ is continuous in a neighborhood of $(-t_0,y_0)$, and therefore problem (\ref{ivpdos}) has some solution $\phi$ defined for $t \in [-t_0,-t_0+\varepsilon_2]$  ($\varepsilon_2>0$). Hence $ \varphi(t)=\phi(-t)$ solves (\ref{ivp}) in the interval $[t_0-\varepsilon_2,t_0]$. Now it suffices to use any solution defined on the left of $t_0$ and any solution defined on the right to have a solution defined in a neighborhood of $t_0$. The claim is proven.
\bigbreak

In the sequel, $\| \cdot\|$ denotes the maximum norm in $\R^n$, i.e., for a vector $x=(x_1,x_2,\dots,x_n) \in \R^n$ we define $$\|x\|=\max_{1 \le j \le n}|x_j|.$$
 
\noindent
{\bf The proof.} As usual, we start by fixing some constants $a>0$ and $b>0$ such that the function $f$ is defined and continuous in the $(n+1)$--dimensional interval $[t_0,t_0+a] \times \{y \in \R^n \, : \, \|y-y_0\| \le b\}$, which is a compact subset of $\R^{n+1}$. Hence there exists $L>0$ such that
$$\|f(t,y)\|\le L \quad \mbox{whenever $0\le t-t_0 \le a$ and $\|y-y_0\|\le b$.}$$

We now define the real interval $I=[t_0,t_0+c]$ with length
 $$c=\min\{a,b/L\},$$ 
and we consider the set ${\cal A}$ of all functions $\gamma:I \longrightarrow \R^n$ such that $\gamma(t_0)=y_0$ and which satisfy a Lipschitz condition with constant $L$, i.e.,
$$\|\gamma(t)-\gamma(s)\| \le L |t-s| \quad \mbox{for all $s,t \in I$.}$$

The previous choice of the constant $c$ guarantees that every function $\gamma \in {\cal A}$ satisfies $\|\gamma(t)-y_0\|\le b$ for all $t \in I$. Hence, for every  $\gamma \in {\cal A}$ the composition $t \in I \longmapsto f(t,\gamma(t))$ is well--defined, bounded (by $L$) and continuous in $I$. It is therefore possible to construct a mapping $$F:{\cal A}\longrightarrow [0,+\infty)$$ as follows: for each function $\gamma \in {\cal A}$ we define
the non--negative number
$$F(\gamma)=\max_{t \in I}\left\|\gamma(t)-y_0-\int_{t_0}^t{f(s,\gamma(s)) \, ds}\right\|,$$
which is a sort of measure of how far the function $\gamma$ is from being a solution. In fact,
the Fundamental Theorem of Calculus ensures that if $F(\gamma)=0$ for some $\gamma \in {\cal A}$ then $\gamma$ is a solution of the initial value problem (\ref{ivp}) in the interval $I$ (the converse is also true, but we do not need it for this proof.)

It is easy to check that the mapping $F$ is continuous in ${\cal A}$ (equipped with the topology of the uniform convergence in $I$). In turn, by the Arzel\`a-- Ascoli Theorem, the domain ${\cal A}$ is  compact. Hence $F$ attains a minimum at some $\varphi \in {\cal A}$. 

To show that $F(\varphi)=0$ (which implies that $\varphi$ is a solution) it suffices to prove that $F$ assumes arbitrarily small positive values in ${\cal A}$. To do so, we follow Tonelli, and for $k \in \N$, $k \ge 2$, we consider the approximate problem\footnote{The differential equations in the approximate problems belong to the class of differential equations {\it with delay}. See, for instance, \cite{sm}.}
$$\left\{
\begin{array}{l}
\mbox{$y(t)=y_0$ for all $t \in [t_0,t_0+c/k]$,}\\
\\
y'(t)=f(t-c/k,y(t-c/k)) \quad \mbox{for all $t \in (t_0+c/k, t_0+c]$.}
\end{array}
\right.$$
This problem has a unique solution $\gamma_k \in {\cal A}$ which we can integrate. Indeed, using induction with the subdivision $t_0,t_0+c/k,t_0+2c/k, \dots,t_0+c$, and changing variables in the corresponding integrals, one can prove that the unique solution satisfies $\|\gamma_k(t)-y_0\| \le b$ ($t \in I$) and
$$\gamma_k(t)=\left\{
\begin{array}{cl}
y_0 & \mbox{for all $t \in [t_0,t_0+c/k]$,}\\
\\
y_0+\int_{t_0}^{t-c/k}{f(s,\gamma_k(s))\, ds} &  \mbox{for all $t \in (t_0+c/k, t_0+c]$.}
\end{array}
\right.$$

 Therefore, for $t \in [t_0,t_0+c/k]$ we have
$$\left\|\gamma_k(t)-y_0-\int_{t_0}^{t}{f(s,\gamma_k(s)) \, ds}\right\|=\left\|\int_{t_0}^{t}{f(s,y_0) \, ds}\right\| \le \dfrac{Lc}{k},$$
and for $t \in (t_0+c/k,t_0+c]$ we have
$$\left\|\gamma_k(t)-y_0-\int_{t_0}^{t}{f(s,\gamma_k(s)) \, ds}\right\|=\left\|\int_{t-c/k}^{t}{f(s,\gamma_k(s)) \, ds}\right\| \le \dfrac{ Lc}{k}.$$
Hence
$$0 \le F(\varphi) \le F(\gamma_k) \le \dfrac{ L c}{k},$$
thus proving that $F(\varphi)=0$ because $k$ can be chosen as big as we wish. \qed
 
The sequence $\{\gamma_k\}_{k \in \Nn}$ constructed in the proof is often referred to as a Tonelli sequence. It is possible to use some other types of minimizing sequences in our proof. In particular, the usual Euler--Cauchy polygons are adequate for the proof too. 

\begin{remark}
\label{rem2}
The proof gives us more information than that collected in the statement, which we have preferred to keep in that form for simplicity and clarity.  

For completeness and for later purposes, we now point out some consequences. In the conditions of Peano's Theorem, and with the notation introduced in the proof, the following results hold:
\begin{enumerate}
\item If $f:[t_0,t_0+a] \times \{y \in \R^n \, : \, \|y-y_0\| \le b\} \to \R^n$ is continuous, bounded by $L>0$ on its domain, and $a \le b/L$, then problem (\ref{ivp}) has at least one solution defined on the whole interval $[t_0,t_0+a]$.
\item A specially important consequence of the previous result arises when $f:[t_0,t_0+a]\times \R^n \longrightarrow \R^n$ is continuous and bounded. In that case we can guarantee that the initial value problem (\ref{ivp}) has at least one solution defined on $[t_0,t_0+a]$.
\end{enumerate}
\end{remark}

Peano's Theorem allows the existence of infinitely many solutions. We owe the following example precisely to Peano.

\bigbreak

\noindent
{\bf Peano's example of a problem with infinitely many solutions.} The scalar problem
\begin{equation}
\label{expe}
y'=3y^{2/3}, \quad y(0)=0,
\end{equation}
has infinitely many solutions. Indeed, one can easily check that $\varphi(t)=0$ and $\psi(t)=t^3$ ($t \in\R$) are solutions. Now the remaining solutions are given by, let us say, a combination of those two. Specifically, if $t_1 \le 0 \le t_2$ then the function
$$\phi(t)=\left\{
\begin{array}{cl}
(t-t_1)^3, & \mbox{if $t < t_1$,}\\
0, & \mbox{if $t_1 \le t \le t_2$,} \\
(t-t_2)^3, & \mbox{if $t >t_2$,}
\end{array}
\right.$$
is a solution of (\ref{expe}). The converse is true as well: every solution of (\ref{expe}) is one of those indicated above.
\section{A finer result in dimension one}

In this section we are concerned with a not so well--known version of Peano's Theorem in dimension $n=1$ which ensures the existence of the least and the greatest solutions to (\ref{ivp}). This result goes back precisely to Peano \cite{peano}, where the greatest solution defined on the right of $t_0$ was obtained as the infimum of all {\it strict upper solutions}. A strict upper solution to the problem (\ref{ivp}) on some interval $I$ is, roughly speaking, some function $\beta=\beta(t)$ satisfying $\beta'(t)>f(t,\beta(t))$ for all $t \in I$, and $\beta(t_0)\ge y_0$. 
Peano also showed in \cite{peano} that the least solution is the supremum of all strict lower solutions, which we define by reversing all the inequalities in the definition of strict upper solution.

\bigbreak

Notice, for instance, that Peano's example (\ref{expe}) has infinitely many solutions, the least one being
$$\mbox{$\varphi_*(t)=t^3$ for $t<0$,} \quad \mbox{$\varphi_*(t)=0$ for $t \ge 0$,}$$
and the greatest solution being
$$\mbox{$\varphi^*(t)=0$ for $t<0$,} \quad \mbox{$\varphi^*(t)=t^3$ for $t \ge 0.$}$$

\bigbreak

Here we present a very easy proof of the existence of the least and the greatest solutions which does not lean on lower/upper solutions (as in \cite{dv, goo, peano, perron}) or on special sequences of approximate solutions (as in \cite{dv, jwal, wwal}). Basically, we obtain the greatest solution as the solution having the greatest integral. This idea works in other settings, see \cite{fp}.

\begin{theorem}
\label{t2} (Second version of Peano's Theorem)
Consider problem (\ref{ivp}) in dimension $n=1$, and assume that there exist constants $a, \, b, \, L \in (0,+\infty)$ such that the function 
$$f:[t_0,t_0+a] \times [y_0-b,y_0+b] \longrightarrow \R$$ is continuous and  $|f(t,y)| \le L$ for all $(t,y) \in [t_0,t_0+a] \times [y_0-b,y_0+b]$.

Then there exist solutions of (\ref{ivp}) $\varphi_*,\varphi^*:I=[t_0,t_0+c]  \longrightarrow \R$, where $c=\min\{a,b/L\}$, such that every solution of (\ref{ivp}) $\varphi:I \longrightarrow \R$ satisfies
$$\varphi_*(t) \le \varphi(t) \le \varphi^*(t) \quad \mbox{for all $t \in I$.}$$

\end{theorem}

\noindent
{\bf Proof.} Let us consider the set of functions ${\cal A}$ introduced in the proof of Peano's Theorem in Section 2 and adapted to dimension $n=1$, i.e., the set of all functions $\gamma:I=[t_0,t_0+c]\longrightarrow \R$ such that $\gamma(t_0)=y_0$ and  
$$|\gamma(t)-\gamma(s)| \le L |t-s| \quad \mbox{for all $s,t \in I$.}$$
Let ${\cal S}$ denote the set of solutions of (\ref{ivp}) defined on $I$. Our first version of Peano's Theorem ensures that ${\cal S}$ is not an empty set. Moreover standard arguments show that ${\cal S} \subset {\cal A}$ and that ${\cal S}$ is a compact subset of ${\cal C}(I)$. Hence the continuous mapping
$${\cal I}: \varphi \in {\cal S} \longmapsto {\cal I}(\varphi)=\int_{t_0}^{t_0+c}{\varphi(s) \, ds}$$
attains a maximum at some $\varphi^* \in {\cal S}$. 

Let us show that $\varphi^*$ is the greatest solution of (\ref{ivp}) on the interval $I$. Reasoning by contradiction, assume that we have a solution $\varphi:I \longrightarrow \R$ such that $\varphi(t_1)> \varphi^*(t_1)$ for some $t_1 \in (t_0,t_0+L)$. Since $\varphi(t_0)=y_0=\varphi^*(t_0)$ we can find
$t_2 \in [t_0,t_1)$ such that $\varphi(t_2)=\varphi^*(t_2)$ and $\varphi> \varphi^*$ on $(t_2,t_1]$. We now have two possibilites on the right of $t_1$: either $\varphi > \varphi^*$ on $(t_2,t_0+L)$, or there exists $t_3 \in (t_1,t_0+L)$ such that $\varphi> \varphi^*$ on $(t_2,t_3)$ and $\varphi(t_3)=\varphi^*(t_3)$. Let us assume the latter (the proof is similar in the other situation) and consider the function 
$$\varphi_1: t \in I \longmapsto \varphi_1(t)=\left\{ 
\begin{array}{cl}
\varphi(t), & \mbox{if $t \in [t_2,t_3]$,}\\
\\
\varphi^*(t), & \mbox{otherwise.}
\end{array}
\right.$$

Elementary arguments with side derivatives show that $\varphi_1 \in {\cal S}$. Moreover $\varphi^* \le \varphi_1$ in $I$, with strict inequality in a subinterval, hence
$${\cal I}(\varphi^*) < {\cal I}(\varphi_1),$$
but this is a contradiction with the choice of $\varphi^*.$

Similarly, one can prove that ${\cal I}$ attains a minimum at certain $\varphi_* \in {\cal S}$, and that $\varphi_*$ is the least element in ${\cal S}$.
\qed

%
%

\medbreak

Can Theorem \ref{t2} be adapted to systems? Yes, it can, but more than continuity must be required for the function $f$, as we will specify below. The need of some extra conditions is easily justified with examples of the following type.
\begin{example}
\label{ex1}
Consider the system
$$\left\{
\begin{array}{ll}y_1'=3y_1^{2/3}, & y_1(0)=0,\\
\\
y_2'=-y_1, &   y_2(0)=0.
\end{array}
\right.$$
The first problem can be solved independently, and it has infinitely many solutions (this is Peano's example again). Notice that the greater the solution we choose for $y_1$ is then the smaller the corresponding $y_2$ becomes on the right of $t_0=0$. Therefore this system does not have a solution which is greater than the other ones {\it in both components.}
\end{example}

Notice however that the system in Example \ref{ex1} has a solution whose first component is greater than the first component of any other solution, and the same is true replacing ``greater" by ``smaller" or ``first component" by ``second component". This observation leads us naturally to the following question: in the conditions of Peano's Existence Theorem we fix a component $i \in \{1,2,\dots,n\}$, can we ensure the existence of a solution with the greatest $i$-th component? The following example answers this question on the negative.

\begin{example}
Let $\phi:\R \longrightarrow \R$ be a continuously differentiable function such that $\phi(0)=0$ and $\phi$ assumes both negative and positive values in every neighborhood of $0$ (hence $\phi'(0)=0$)\footnote{A paradigm is $\phi(x)=x^3 \sin(1/x)$ for $x \neq 0$, and $\phi(0)=0$.}.

The idea is to construct a system whose solutions have a component which is a translation of $\phi$, and then those specific components cannot be compared. To do so it suffices to consider the two dimensional system
$$\left\{
\begin{array}{ll}y_1'=3y_1^{2/3}, & y_1(0)=0,\\
\\
y_2'=\phi' \left(y_1^{1/3} \right), &   y_2(0)=0.
\end{array}
\right.$$
 
Let $\varepsilon >0$ be fixed; we are going to prove that there is not a solution of the system whose second component is greater than the second component of any other solution on the whole interval $[0,\varepsilon]$.

First, note that we can compute all the solutions. For each $a \in [0,\varepsilon]$ we have a solution $\varphi(t)=(\varphi_1(t),\varphi_2(t))$ given by
$$\mbox{$\varphi_1=0$ on $[0,a]$} \quad \mbox{and} \quad \mbox{$\varphi_1(t)=(t-a)^3$ for $t  \in [a,\varepsilon]$,}$$
and then it suffices to integrate the second component to obtain $\varphi_2=0$ on $[0,a]$, and
$$\varphi_2(t)=\int_b^t{\phi'(s-a) \, ds}=\phi(t-a) \quad \mbox{for $t \in [a,\varepsilon]$.}$$
Conversely, every solution of the system is given by the previous expression for an adequate value of $a \in [0,\varepsilon]$. 

Now let us consider two arbitrary solutions of the system. They are given by the above formulas for some corresponding values $a=b$ and $a=b'$, with $0 \le b<b'\le \varepsilon$, and then their respective second components cannot be compared in the subinterval $(b,b')$.
 \end{example}

The previous example still has a solution with the greatest first component, but, in the author's opinion, this is just a consequence of the fact that the first equation in the system is uncoupled and we can solve it independently (in particular, Theorem \ref{t2} applies). However this remark raises the open problem of finding a two dimensional system which has neither a solution with the greatest first component nor a solution with the greatest second component.
 
 \bigbreak

A multidimensional version of Theorem \ref{t2} is valid if the nonlinear part
$$f(t,y)=(f_1(t,y),f_2(t,y),\dots, f_n(t,y))$$
is {\it quasimonotone nondecreasing}, i.e., if for each component $i \in \{1,2,\dots,n\}$ the relations $y_j \le \bar y_j$, $j \neq i$, imply
$$f_i(t,(y_1,\dots,y_{i-1},y,y_i,\dots,y_n)) \le f_i(t,(\bar y_1,\dots,\bar y_{i-1},y,\bar y_i,\dots,\bar y_n)).$$
The reader is referred to \cite{bp, bs, has, hl, wal, wwal} and references therein for more information on quasimonotone systems.

%
%
%
%
%

\section{The power of lower and upper solutions: Existence for nonlocal problems}
The real power of lower and upper solutions reveals when we want to guarantee the existence of solution to (\ref{ivp}) on a given interval, and not merely in an unknown (possibly very small) neighborhood of $t_0$. 

Let $a>0$ be fixed, let $f:[t_0,t_0+a]\times \R \longrightarrow \R$ be continuous, and consider the nonlocal problem
\begin{equation}
\label{ivp2}
y'=f(t,y) \quad \mbox{for all $t \in I=[t_0,t_0+a]$,}Ê\quad y(t_0)=y_0.
\end{equation}

In this section we follow Goodman \cite{goo}, who considered the Carath\'eodory version of (\ref{ivp2}) and proved that the greatest solution is the supremum of all (nonstrict) lower solutions, while the least solution is the infimum of all upper solutions.\footnote{Note that we can find in the literature some other denominations for upper (lower) solutions, such as upper (lower) functions, or superfunctions (subfunctions).} 
In this paper, by upper solution we mean a function $\beta:I \longrightarrow \R$ which is continuously differentiable in the interval $I$, $\beta(t_0) \ge y_0$, and $\beta'(t) \ge f(t, \beta(t))$ for all $t \in I$. We define a lower solution in an analogous way reversing the corresponding inequalities.

\bigbreak

Our first result is quite simple and well--known. Indeed, it has a standard concise proof as a corollary of Peano's Existence Theorem (one has to use the second observation in Remark \ref{rem2}). However it can also be proven independently, as we are going to show, and then Peano's Existence Theorem in dimension $n=1$ will follow as a corollary (see Remark \ref{remp}). To sum up, the following proof is another new proof to Peano's Theorem in dimension $n=1$ which readers might want to compare with those in \cite{dv, gar, jwal, wwal}. The main difference with respect to the proofs in \cite{dv, gar, jwal, wwal} is that we need to produce approximate solutions between given lower and upper solutions, which we do with the aid of Lemma \ref{lema}.

\begin{theorem}
\label{tss1} (Third version of Peano's Theorem)
Suppose that problem (\ref{ivp2}) has a lower solution $\alpha$ and an upper solution $\beta$ such that $\alpha(t) \le \beta(t)$ for all $t \in I$.  

Then problem (\ref{ivp2}) has at least one solution $\varphi:I \longrightarrow \R$ such that  
$$ \alpha(t) \le \varphi(t) \le \beta(t) \quad \mbox{for all $t \in I$.}$$
\end{theorem}

\noindent
{\bf Proof.} Let $L>0$ be fixed so that
$$|f(t,y)| \le L \quad \mbox{for all $(t,y) \in I \times \R^2$ such that $\alpha(t) \le y \le \beta(t)$,}$$
and $\max \{|\alpha'(t)|, |\beta'(t)| \} \le L$ for all $t \in I$. Let us define the set ${\cal A}$ of all continuous functions $\gamma:I \longrightarrow \R$ such that $\alpha \le \gamma \le \beta$ on $I$ and
$$|\gamma(t) -\gamma(s)| \le L |t-s| \quad \mbox{for all $s,t \in I$.}$$
The choice of $L$ ensures that $\alpha, \, \beta \in {\cal A}$, so ${\cal A}$ is not empty. Moreover, the set ${\cal A}$ is a connected subset of ${\cal C}(I)$ (convex, actually), and the Arzel\`a--Ascoli Theorem implies that ${\cal A}$ is compact. Hence the mapping  defined by
$$F(\gamma)=\max_{t \in I}\left|\gamma(t)-y_0-\int_{t_0}^t{f(s,\gamma(s)) \,ds} \right| \quad \mbox{for each $\gamma \in {\cal A},$}$$
attains a minimum at some $\varphi \in {\cal A}$. We are going to prove that $F(\varphi)=0$, thus showing that $\varphi$ is a solution in the conditions of the statement. To do it, we are going to prove that  we can find functions $\gamma \in {\cal A}$ such that $F(\gamma)$ is as small as we wish. The construction of such functions leans on the following lemma.

\begin{lemma}
\label{lema}
For all $t_1, t_2 \in I$ such that $t_1 < t_2$, and all $y_1 \in [\alpha(t_1),\beta(t_1)]$, there exists $\gamma \in {\cal A}$ such that
$$\gamma(t_1)=y_1 \quad \mbox{and} \quad \gamma(t_2)=y_1+\int_{t_1}^{t_2}{f(s,\gamma(s)) \, ds}.$$
\end{lemma}

\noindent
{\bf Proof of Lemma \ref{lema}.} Let $t_1,t_2$ and $y_1$ be as in the statement. We define a set of functions ${\cal A}_1=\{\gamma \in {\cal A} \, : \, \gamma(t_1)=y_1 \}$.
 The set ${\cal A}_1$ is not empty: an adequate convex linear combination of $\alpha$ and $\beta$ assumes the value $y_1$ at $t_1$, and then it belongs to ${\cal A}_1$.

Let us consider the mapping $G: {\cal A}_1 \longrightarrow \R$, defined for each $\gamma \in {\cal A}_1$ as
$$G(\gamma)=\gamma(t_2)-y_1-\int_{t_1}^{t_2}{f(s,\gamma(s)) \, ds}.$$
To finish the proof it suffices to show that $G(\gamma)=0$ for some $\gamma \in {\cal A}_1$.

The mapping $G$ is continuous in ${\cal A}_1$, which is connected, hence $G({\cal A}_1)$ is a connected subset of the reals, i.e., an interval\footnote{This is not true in dimension $n>1$, so this approach does not work in that case.}. Therefore to ensure the existence of some $\gamma \in {\cal A}_1$ such that $G(\gamma)=0$ it suffices to prove the existence of functions $\tilde \alpha$ and $\tilde \beta$ in ${\cal A}_1$ such that $G(\tilde \alpha) \le 0 \le G(\tilde \beta)$. 

Next we show how to construct one such $\tilde \alpha$ from $\alpha$ (the construction of $\tilde \beta$ from $\beta$ is analogous and we omit it). If $\alpha(t_1)=y_1$ we simply take $\tilde \alpha=\alpha$. If, on the other hand, $\alpha(t_1) < y_1$ then we define $\tilde \alpha$ in ``three (or two) pieces": first, we define $\tilde \alpha$ on $[t_0,t_1]$ as an adequate convex linear combination of $\alpha$ and $\beta$ to have $\tilde \alpha(t_1)=y_1$; second, we define $\tilde \alpha$ on the right of $t_1$ as the function whose graph is the line with slope $-L$ starting at the point $(t_1,y_1)$ and on the interval $[t_1,t_3]$, where $t_3$ is the first point in the interval $(t_1,t_0+a)$ such that the line intersects with the graph of $\alpha$, and finally we continue $\tilde \alpha=\alpha$ on $[t_3,t_0+a]$. If no such $t_3$ exists, then $\tilde \alpha$ is simply the line with slope $-L$ on the whole interval $[t_1,t_0+a]$. Verifying that $\tilde \alpha \in {\cal A}_1$ and that $G(\tilde \alpha) \le 0$ is just routine. The proof of Lemma \ref{lema} is complete. \qed

\bigbreak

Now we carry on with the final part of the proof of Theorem \ref{tss1}. Let $\varepsilon>0$ be fixed and consider a partition of the interval $[t_0,t_0+a]$, say $t_0,t_1, \dots, t_k=t_0+a$ ($k \in \N$), such that
$$0<t_j-t_{j-1}<\dfrac{\varepsilon}{2L} \quad \mbox{for all $j \in \{1,2,\dots,k\}$.}$$

 Lemma \ref{lema} guarantees that we can construct (piece by piece) a function $\gamma_{\varepsilon} \in {\cal A}$ such that $\gamma_{\varepsilon}(t_0)=y_0$ and
$$\gamma_{\varepsilon}(t_j)=\gamma_{\varepsilon}(t_{j-1})+\int_{t_{j-1}}^{t_j}{f(s,\gamma_{\varepsilon}(s)) \, ds} \quad \mbox{for all $j \in \{1,2,\dots, k\}$.}$$

Now for each $t \in (t_0,t_0+a]$ there is a unique $j \in \{1,2,\dots, k\}$ such that $t \in (t_{j-1},t_j]$ and then
\begin{align*}
\left| \gamma_{\varepsilon}(t)-y_0-\int_{t_0}^t{f(s,\gamma_{\varepsilon}(s)) \, ds} \right|&=\left| \gamma_{\varepsilon}(t)-\gamma_{\varepsilon}(t_{j-1})-\int_{t_{j-1}}^{t}{f(s,\gamma_{\varepsilon}(s)) \, ds} \right|\\
& \le |\gamma_{\varepsilon}(t)-\gamma_{\varepsilon}(t_{j-1})|+\int_{t_{j-1}}^{t}{|f(s,\gamma_{\varepsilon}(s))| \, ds}\\
& \le 2L|t-t_{j-1}| < \varepsilon.
\end{align*}
Hence $0 \le F(\varphi) \le F(\gamma_{\varepsilon}) < \varepsilon$, which implies that $F(\varphi)=0$ because $\varepsilon>0$ was arbitrarily chosen. \qed

\begin{remark}
\label{remp}
Notice that Peano's Existence Theorem in dimension $n=1$ is really a consequence of Theorem \ref{tss1}. To see it consider problem (\ref{ivp}) in dimension $n=1$ and let $b>0$, $L>0$ and $c>0$ be as in the proof of Peano's Theorem in Section 2. Now define a new function
$$\tilde f(t,y)=f(t,\max\{ y_0-b,\min\{y,y_0+b\}\}) \quad \mbox{for all $(t,y) \in [t_0,t_0+c]\times \R$,}$$
which is continuous and bounded by $L$ on the whole of $[t_0,t_0+c]  \times \R$.

Obviously, the functions $\alpha(t)=y_0-L(t-t_0)$ and $\beta(t)=y_0+L(t-t_0)$ are, respectively, a lower and an upper solution to
\begin{equation}
\nonumber
y'=\tilde f(t,y), \, \, t \in [t_0,t_0+c], \quad  y(t_0)=y_0,
\end{equation}
which then has at least one solution $\varphi \in [\alpha,\beta]$, by virtue of Theorem \ref{tss1}. Since $|\varphi'(t)| \le L$ for all $t  \in [t_0,t_0+c]$ we have $|\varphi(t)-y_0| \le b$ for all $t \in [t_0,t_0+c]$, and then the definition of $\tilde f$ implies that $\varphi$ solves (\ref{ivp}).
\end{remark}

The existence of both the lower and the upper solutions in Theorem \ref{tss1} is essential. Our next example shows that we cannot expect to have solutions {\it for a nonlocal problem} if we only have a lower (or an upper) solution.

\begin{example}
The function $\alpha(t)=0$ for all $t \in [0,\pi]$ is a lower solution to the initial value problem
$$y'=1+y^2, \quad y(0)=0,$$
which has no solution defined on $[0,\pi]$ (Its unique solution is $\varphi(t)=\mbox{\rm tan} \, t$ for all $t \in (-\pi/2, \pi/2)$).
\end{example} 

Theorem \ref{tss1} does not guarantee that every solution of (\ref{ivp2}) is located between the lower and the upper solutions. To see it, simply note that solutions are lower and upper solutions at the same time, so the zero function is both lower and upper solution for Peano's example (\ref{expe}) which has many solutions above the zero function.  

However it is true that, in the conditions of Theorem \ref{tss1}, we have a greatest and a least solution between $\alpha$ and $\beta$, and we have their respective Goodman's characterizations in terms of lower and upper solutions, see \cite{goo}.

\begin{corollary}
\label{t4}

Suppose that problem (\ref{ivp2}) has a lower solution $\alpha$ and an upper solution $\beta$ such that $\alpha(t) \le \beta(t)$ for all $t \in I$ and let
$$[\alpha, \beta]=\{\gamma \in {\cal C}(I) \, :\, \mbox{$\alpha \le \gamma \le \beta$ on $I$} \}.$$

Then problem (\ref{ivp2}) has solutions $\varphi_*,\varphi^* \in [\alpha,\beta]$ such that every solution of (\ref{ivp}) $\varphi \in [\alpha,\beta]$ satisfies
$$ \varphi_*(t) \le \varphi(t) \le \varphi^*(t) \quad \mbox{for all $t \in I$.}$$

Moreover, the least solution of (\ref{ivp2}) in $[\alpha,\beta]$ satisfies
\begin{equation}
\label{le}
\varphi_*(t)=\min \{ \gamma(t) \, : \, \mbox{$\gamma \in [\alpha,\beta]$, $\gamma$ upper solution of (\ref{ivp2})}\} \quad (t \in I),
\end{equation}
and the greatest solution of (\ref{ivp}) in $[\alpha,\beta]$ satisfies
\begin{equation}
\label{le2}
\varphi^*(t)=\max \{ \gamma(t) \, :	\, \mbox{$\gamma \in [\alpha,\beta]$, $\gamma$ lower solution of (\ref{ivp2})}\} \quad (t \in I).
\end{equation}
\end{corollary}

\noindent
{\bf Proof.} Let $L>0$ be fixed so that
$$|f(t,y)| \le L \quad \mbox{for all $(t,y) \in I \times \R^2$ such that $\alpha(t) \le y \le \beta(t)$.}$$
 
The set of solutions of (\ref{ivp2}) in $[\alpha,\beta]$ is not empty by virtue of Theorem \ref{tss1}. Moreover, it is a compact subset of ${\cal C}(I)$, because it is closed, bounded and every one of its elements satisfies a Lipschitz condition with constant $L$. Hence, a similar argument to that in the proof of Theorem \ref{t2} guarantees that there is a solution $\varphi^* \in [\alpha,\beta]$ which is greater than any other solution $\varphi \in [\alpha,\beta]$.

Let us prove that $\varphi^*$ satisfies (\ref{le2}). Notice first that if $\gamma \in [\alpha,\beta]$ is a lower solution to (\ref{ivp2}) then Theorem \ref{tss1} guarantees that there is some solution $\varphi_{\gamma} \in [\gamma,\beta]$. Hence, for all $t \in I$ we have
\begin{align*}
\sup\{\gamma(t) \, :\, \mbox{$\gamma \in [\alpha,\beta]$, $\gamma$ lower solution}  \} &\le \max\{\varphi(t) \, :\, \mbox{$\varphi \in [\alpha,\beta]$, $\varphi$ solution}
 \}\\
&=\varphi^*(t),
\end{align*}
and then (\ref{le2}) obtains because $\varphi^*$ is a lower solution of (\ref{ivp2}) in $[\alpha,\beta]$.

The proof of the existence of $\varphi_*$ and the proof of (\ref{le}) are similar.
\qed

\begin{remark}
Corollary \ref{t4} can be proven directly via Perron's method, see \cite{dv, goo, po2, perron}. This means starting at (\ref{le2}) as a definition and then showing that $\varphi^*$ is a solution in the conditions of the statement. Perron's method involves careful work with sets of one--sided differentiable functions, which we avoid. Corollary \ref{t4} can also be proven easily from Theorem \ref{t2}.
\end{remark}
 
Finally we deduce from Corollary \ref{t4} the Peano's characterizations of the least and the greatest solutions in terms of strict lower and upper solutions. In doing so we are finally proving the ``real" Peano's Theorem, because our next result is the closest to the one proven in \cite{peano}.

We say that $\alpha:I \longrightarrow \R$ is a strict lower solution of (\ref{ivp2}) if it is continuously differentiable in $I$, $\alpha(t_0)\le y_0$, and $\alpha'(t)<f(t,\alpha(t))$ for all $t \in I$. A strict upper solution is defined analogously by reversing the relevant inequalities. Notice that if $\alpha$ is a strict lower solution and $\beta$ is a strict upper solution, then $\alpha < \beta$ on $(t_0,t_0+a]$, for otherwise we could find $t_1 \in (t_0,t_0+a]$ such that $\alpha< \beta$ in $(t_0,t_1)$ and $\alpha(t_1)=\beta(t_1)$, but then we would have
$$\beta'(t_1)>f(t_1,\beta(t_1))=f(t_1,\alpha(t_1)) > \alpha'(t_1) \ge \beta'(t_1), \quad \mbox{a contradiction.}$$

\begin{theorem}
\label{t5} (Fourth version of Peano's Theorem)
If problem (\ref{ivp2}) has a strict lower solution $\alpha$ and a strict upper solution $\beta$ then the following results hold:
\begin{enumerate}
\item Problem (\ref{ivp2}) has at least one solution;
\item If $\varphi$ solves (\ref{ivp2}) just in some domain $[t_0,t_0+\varepsilon]$, with $\varepsilon \in (0,a)$, then it can be extended as a solution of (\ref{ivp2}) to the whole interval $[t_0,t_0+a]$ and then it satisfies
\begin{equation}
\label{dee}
\mbox{$\alpha(t) < \varphi(t) < \beta(t)$ for all $t \in (t_0,t_0+a]$;}
\end{equation} 
\item Problem (\ref{ivp2}) has the least solution $\varphi_*:I  \longrightarrow \R$ and the greatest solution $\varphi^*:I \longrightarrow \R$, which satisfy
\begin{equation}
\label{lepe}
\varphi_*(t)=\sup \{ \gamma(t) \, : \, \mbox{ $\gamma$ strict lower solution of (\ref{ivp2})}\} \quad (t \in I),
\end{equation}
and  
\begin{equation}
\label{le2pe}
\varphi^*(t)=\inf \{ \gamma(t) \, :	\, \mbox{ $\gamma$ strict upper solution of (\ref{ivp2})}\} \quad (t \in I).
\end{equation}
\end{enumerate}
\end{theorem}

\noindent
{\bf Proof.} Corollary \ref{t4} guarantees that (\ref{ivp2}) has the least and the greatest solutions in $[\alpha,\beta]=\{\gamma \in {\cal C}(I) \, : \, \alpha \le \gamma \le \beta\}$, which we denote, respectively, by $\varphi_*$ and $\varphi^*$.
 
Let us prove that solutions of (\ref{ivp2}) satisfy (\ref{dee}) which, in particular, ensures that all of them belong to $[\alpha,\beta]$. Reasoning by contradiction, assume that a certain solution $\varphi:[t_0,t_0+\varepsilon] \longrightarrow \R$ ($\varepsilon \in (0,a)$) satisfies $\varphi(t_1) \ge \beta(t_1)$ for some $t_1 \in (t_0,t_0+\varepsilon)$. The initial conditions ensure that we can find some $t_2 \in (t_0,t_1)$ such that $\varphi<\beta$ in $(t_0,t_2)$ and $\varphi(t_2)=\beta(t_2)$, but then we have
$$\beta'(t_2)>f(t_2,\beta(t_2))=f(t_2,\varphi(t_2))=\varphi'(t_2) \ge \beta'(t_2),$$
a contradiction. This proves that every solution is smaller than $\beta$ on its domain, and a similar argument shows that every solution is greater than $\alpha$. Therefore Theorem \ref{tss1} ensures that every solution of (\ref{ivp2}) can be continued to the whole interval $I=[t_0,t_0+a]$ as a solution of (\ref{ivp2}) and between $\alpha$ and $\beta$. Hence (\ref{dee}) obtains and, moreover, $\varphi_*$ and $\varphi^*$ are, respectively, the least and the greatest among all the solutions of (\ref{ivp2}).

Next we show that (\ref{le2pe}) is satisfied. The proof of (\ref{lepe}) is similar and we omit it.

The previous arguments still work if we replace $\beta$ by any other strict upper solution $\gamma:I \longrightarrow \R$ (necessarily greater than $\alpha$ on $I$). Hence
$$\varphi^*(t) \le \inf \{ \gamma(t) \, :	\, \mbox{$\gamma$ strict upper solution of (\ref{ivp2})}\} \quad (t \in I).$$
To show that we can replace this inequality by an identity it suffices to construct a decreasing sequence of strict upper solutions which converges to $\varphi^*$. To do it, let $k \in \N$ be sufficiently large so that
$$\beta'(t)>f(t,\beta(t))+1/k \quad \mbox{for all $t \in I$.}$$
Plainly, $\beta$ is an upper solution to the initial value problem
\begin{equation}
\label{pax}
y'=f(t,y)+1/k, \quad y(t_0)=y_0,
\end{equation}
and, in turn, $\varphi^*$ is a lower solution. Hence, for all sufficiently large values of $k$ there exists $\varphi_k$, a solution of (\ref{pax}), between $\varphi^*$ and $\beta$.

Obviously, $\varphi_k$ is a strict upper solution to (\ref{pax}) with $k$ replaced by $k+1$, hence $\varphi_k \ge \varphi_{k+1} \ge \varphi^*$ on $I$. Thus we can define a limit function
\begin{equation}
\label{ga}
\varphi_{\infty}(t)= \lim_{k \to \infty}\varphi_k(t) \ge \varphi^*(t) \quad (t \in I).
\end{equation}
The functions $\varphi_k$ satisfy Lipschitz conditions on $I$ with the same Lipschitz constant, which implies that $\varphi_{\infty}$ is Lipschitz continuous on $I$. Dini's Theorem ensures then that the sequence $\{\varphi_k\}_k$ converges uniformly to $\varphi_{\infty}$ on $I$. Now for all sufficiently large values of $k \in \N$ we have 
$$\varphi_k(t)=y_0+\int_{t_0}^t{f(s,\varphi_k(s)) \, ds}+(t-t_0)/k \quad (t \in I),$$
and then taking limit when $k$ tends to infinity we deduce that $\varphi_{\infty}$ is a solution of (\ref{ivp2}). Hence $\varphi_{\infty} \le \varphi^*$ on $I$, and then (\ref{ga})  yields $\varphi_{\infty}=\varphi^*$ on $I$. The proof of (\ref{le2pe}) is complete. \qed

\section{Concluding remarks}
\begin{enumerate}
\item The assumption $\alpha \le \beta$ on $I$ in Theorem \ref{tss1} can be omitted. Marcelli and Rubbioni proved in \cite{mr} that we have solutions between the minimum of $\alpha$ and $\beta$ and the maximum of them. Furthermore, we do not even need that lower or upper solutions be continuous, see \cite{po}.
\item Theorem \ref{tss1} and Corollary \ref{t4} can be extended to quasimonotone systems, see \cite{bp, bs, hl}. Theorem \ref{tss1} for systems even works when the lower and upper solutions are not ordered, see \cite{mr}, and $f$ may be discontinuous or singular, see \cite{bp}.
\item The lower and upper solutions method is today a most acknowledged effective tool in the analysis of differential equations and, specially, boundary value problems. A detailed account on how far the method has evolved (just for second--order ODEs!) is given in the monograph by De Coster and Habets \cite{ch}. As far as the author is aware, the first use of lower and upper solutions is due Peano in \cite{peano}.

\end{enumerate}

\end{document}